\documentclass{amsproc}
\usepackage{amscd}
\theoremstyle{plain}
\newtheorem{theorem}{Theorem}[section]
\theoremstyle{definition}
\newtheorem{definition}[theorem]{Definition}
\newtheorem{question}[theorem]{Question}

\theoremstyle{remark}
\newtheorem{example}[theorem]{Example}
\newtheorem{remark}[theorem]{Remark}
\newcommand{\comps}{{\mathbb C}}
\newcommand{\reals}{{\mathbb R}}
\newcommand{\A}{{\mathfrak A}}
\newcommand{\cstar}{C^\ast}
\newcommand{\tensor}{\otimes}
\newcommand{\Tensor}{\hat{\otimes}}
\newcommand{\VN}{\operatorname{VN}}
\newcommand{\posints}{{\mathbb N}}
\newcommand{\M}{{\mathfrak M}}
\newcommand{\N}{{\mathfrak N}}
\newcommand{\id}{\operatorname{id}}
\newcommand{\wstar}{W^\ast}
\begin{document}
\title[Abstract harmonic analysis]{Abstract harmonic analysis, homological algebra, and operator spaces}
\author{Volker Runde}
\address{Department of Mathematical and Statistical Sciences, University of Alberta, Edmonton, AB, Canada T6G 2G1}
\email{vrunde@ualberta.ca}
\thanks{Financial support by NSERC under grant no.\ 227043-00 is gratefully acknowledged.}
\keywords{locally compact groups, group algebra, Fourier algebra, Fourier--Stieltjes algebra, Hochschild cohomology, homological algebra, operator spaces}
\subjclass{22D15, 22D25, 43A20, 43A30, 46H20 (primary), 46H25, 46L07, 46M18, 46M20, 47B47, 47L25, 47L50}
\begin{abstract}
In 1972, B.\ E.\ Johnson proved that a locally compact group $G$ is amenable if and only if certain Hochschild cohomology groups of its convolution
algebra $L^1(G)$ vanish. Similarly, $G$ is compact if and only if $L^1(G)$ is biprojective: In each case, a classical property of $G$ corresponds to a cohomological 
propety of $L^1(G)$. Starting with the work of Z.-J. Ruan in 1995, it has become apparent that in the non-commutative setting, i.e. when
dealing with the Fourier algebra $A(G)$ or the Fourier--Stieltjes algebra $B(G)$, the canonical operator space structure of the algebras under consideration
has to be taken into account: In analogy with Johnson's result, Ruan characterized the amenable locally compact groups $G$ through the
vanishing of certain cohomology groups of $A(G)$. In this paper, we give a survey of historical developments, known results, and current open problems.
\end{abstract}
\maketitle
\section{Abstract harmonic analysis,...}
The central objects of interest in abstract harmonic analysis are locally compact groups, i.e.\ groups equipped with a locally compact Hausdorff topology such that multiplication
and inversion are continuous. This includes all discrete groups, but also all Lie groups. There are various function spaces associated with a locally compact group $G$, e.g.\ the space
${\mathcal C}_0(G)$ of all continuous functions on $G$ that vanish at infinity. The dual space of ${\mathcal C}_0(G)$ can be identified with $M(G)$, the space of all regular (complex) Borel measures on $G$. The convolution
product $\ast$ of two measures is defined via
\[
  \langle f, \mu \ast \nu \rangle := \int_G \int_G f(xy) \, d\mu(x) \, \nu(y) \qquad (\mu,\nu \in M(G), \, f \in {\mathcal C}_0(G))
\]
and turns $M(G)$ into a Banach algebra. Moreover, $M(G)$ has an isometric involution given by
\[
  \langle f, \mu^\ast \rangle := \int_G \overline{f(x^{-1})} \, d\mu(x) \qquad (\mu \in M(G), \, f \in {\mathcal C}_0(G)).
\]
\par
The most surprising feature of an object as general as a locally compact group is the existence of (left) {\it Haar measure\/}: a regular Borel measure which is invariant under left translation and unique up to a multiplicative constant.
For example, the Haar measure of a discrete group is simply counting measure, and the Haar measure of $\reals^N$, is $N$-dimensional Lebesgue measure. The space $L^1(G)$ of all integrable functions with respect to Haar measure can be identified
with a closed $^\ast$-ideal of $M(G)$ via the Radon--Nikod\'ym theorem. Both $M(G)$ and $L^1(G)$ are complete invariants for $G$: Whenver $L^1(G_1)$ and $L^1(G_2)$ (or $M(G_1)$ and $M(G_2)$) are isometrically isomorphic, then $G_1$ and $G_2$ are
topologically isomorphic. This means that all information on a locally compact group is already encoded in $L^1(G)$ and $M(G)$. For example, $L^1(G)$ and $M(G)$ are abelian if and only if $G$ is abelian, and $L^1(G)$ has an identity if and only
if $G$ is discrete.
\par
References for abstract harmonic analysis are \cite{Fol}, \cite{HR}, and \cite{Rei}.
\par
The property of locally compact groups we will mostly be concerned in this survey is {\it amenability\/}. A a {\it mean\/} on a locally compact group $G$ is a bounded linear functional $m \!: L^\infty(G) \to \comps$ such that 
$\langle 1, m \rangle = \| m \| =1$. For any function $f$ on $G$ and for any $x \in G$, we write $L_x f$ for the left translate of $f$ by $x$, i.e.\ $(L_x f)(y) := f(xy)$ for $y \in G$.
\begin{definition} \label{amdef}
A locally compact group $G$ is called {\it amenable\/} if there is a (left) translation invariant mean on $G$, i.e.\ a mean $m$ such that
\[
  \langle \phi, m \rangle = \langle L_x \phi, m \rangle \qquad (\phi \in L^\infty(G), \, x \in G).
\]
\end{definition}
\begin{example}
\begin{enumerate}
\item Since the Haar measure of a compact group $G$ is finite, $L^\infty(G) \subset L^1(G)$ holds. Consequently, Haar measure is an invariant mean on $G$.
\item For abelian $G$, the Markov--Kakutani fixed point theorem yields an invariant mean on $G$.
\item The free group in two generators is {\it not\/} amenable (\cite[(0.6) Example]{Pat}).
\end{enumerate}
\end{example}
\par
Moreover, amenability is stable under standard constructions on locally compact groups such as taking subgroups, quotients, extensions, and inductive limits.
\par
Amenable, locally compact groups were first considered by J.\ v.\ Neumann (\cite{vonN}) in the discrete case; he used the term ``Gruppen von endlichem Ma{\ss}''. The adjective amenable for groups satisfying Definition \ref{amdef}
is due to M.\ M.\ Day (\cite{Day}), apparently with a pun in mind: They are a{\it men\/}able because they have an invariant mean, but also since they are particularly pleasant to deal with and thus are truly {\it amenable\/} --- just in the sense 
of that adjective in everyday speech.   
\par
For more on the theory of amenable, locally compact groups, we refer to the monographs \cite{Gre}, \cite{Pat}, and \cite{Pie}.
\section{homological algebra,...}
We will not attempt here to give a survey on a area as vast as homological algebra, but outline only a few, basic cohomological concepts that are relevant in connection with abstract harmonic analysis. For the general theory of
homological algebra, we refer to \cite{CE}, \cite{MacL}, and \cite{Wei}. The first to adapt notions from homological algebra to the functional analytic context was H.\ Kamowitz in \cite{Kam}.
\par
Let $\A$ be a Banach algebra. A Banach $\A$-bimodule is a Banach space $E$ which is also an $\A$-bimodule such that the module actions of $\A$ on $E$ are jointly continuous. A {\it derivation\/} from $\A$ to $E$ is a (bounded) 
linear map $D \!: \A \to E$ satisfying
\[
  D(ab) = a \cdot Db + (Da) \cdot b \qquad (a,b \in \A);
\]
the space of all derivation from $\A$ to $E$ is commonly denoted by ${\mathcal Z}^1(\A,E)$. A derivation $D$ is called {\it inner\/} if there is $x \in E$ such that
\[
  Da = a \cdot x - x \cdot a \qquad (a \in \A).
\]
The symbol for the subspace of ${\mathcal Z}^1(\A,E)$ consisting of the inner derivations is ${\mathcal B}^1(\A,E)$; note that ${\mathcal B}^1(\A,E)$ need not be closed in ${\mathcal Z}^1(\A,E)$.
\begin{definition}
Let $\A$ be a Banach algebra, and let $E$ be a Banach $\A$-bi\-mod\-ule. Then then the {\it first Hochschild cohomology group\/} ${\mathcal H}^1(\A,E)$ of $\A$ with coefficients in $E$ is defined as
\[
  {\mathcal H}^1(\A,E) := {\mathcal Z}^1(\A,E) / {\mathcal B}^1(\A,E).
\]
\end{definition}
\par
The name Hochschild cohomology group is in the honor of G.\ Hochschild who first considered these groups in a purely algebraic context (\cite{Hoch1} and \cite{Hoch2}).
\par
Given a Banach $\A$-bimodule $E$, its dual space $E^\ast$ carries a natural Banach $\A$-bimodule structure via
\[
  \langle x, a \cdot \phi \rangle := \langle x \cdot a, \phi \rangle \quad\text{and}\quad \langle x, \phi \cdot a \rangle := \langle a \cdot x, \phi \rangle
  \qquad (a \in \A, \, \phi \in E^\ast, \, x \in E).
\]
We call such Banach $\A$-bimodules {\it dual\/}.
\par
In his seminal memoir \cite{Joh1}, B.\ E.\ Johnson characterized the amenable locally compact groups $G$ through Hochschild cohomology groups of $L^1(G)$ with coefficients in dual Banach $L^1(G)$-bimodules (\cite[Theorem 2.5]{Joh1}):
\begin{theorem}[B.\ E.\ Johnson] \label{Barry1}
Let $G$ be a locally compact group. Then $G$ is amenable if and only if ${\mathcal H}^1(L^1(G),E^\ast) = \{ 0 \}$ for each Banach $L^1(G)$-bimodule $E$.
\end{theorem}
\par
The relevance of Theorem \ref{Barry1} is twofold: First of all, homological algebra is a large and powerful toolkit --- the fact that a certain property is cohomological in nature allows to apply its tools, which then yield
further insights. Secondly, the cohomological triviality condition in Theorem \ref{Barry1} makes sense for {\it every\/} Banach algebra. This motivates the following definition from \cite{Joh1}:
\begin{definition} \label{bamdef}
A Banach algebra $\A$ is called {\it amenable\/} if ${\mathcal H}^1(\A,E^\ast) = \{ 0 \}$ for each Banach $\A$-bimodule $E$.
\end{definition}
\par
Given a new definition, the question of how significant it is arises naturally. Without going into the details and even without defining what a nuclear $\cstar$-algebra is, we would like to only mention the following very deep
result which is very much a collective accomplishment of many mathematicians, among them A.\ Connes, M.\ D.\ Choi, E.\ G.\ Effros, U.\ Haagerup, E.\ C.\ Lance, and S.\ Wassermann:
\begin{theorem}
A $\cstar$-algebra is amenable if and only if it is nuclear.
\end{theorem}
\par
For a relatively self-contained exposition of the proof, see \cite[Chapter 6]{Run}.
\par
Of course, Definition \ref{bamdef} allows for modifications by replacing the class of all dual Banach $\A$-bimodules by any other class. In \cite{BCD}, W.\ G.\ Bade, P.\ C.\ Curtis, Jr., and H.\ G.\ Dales  called a commutative 
Banach algebra $\A$ {\it weakly amenable\/} if and only if ${\mathcal H}^1(\A,E) = \{ 0 \}$ for every symmetric Banach $\A$-bimodule $E$, i.e.\ satisfying
\[
  a \cdot x = x \cdot a \qquad (a \in \A, \, x \in E).
\]
This definition is of little use for non-commutative $\A$. For commutative $\A$, weak amenability, however, is equivalent to ${\mathcal H}^1(\A,\A^\ast) = \{ 0 \}$ (\cite[Theorem 1.5]{BCD}), and in \cite{Joh2}, Johnson suggested that this should be
used to define weak amenability for arbitrary $\A$:
\begin{definition} \label{wamdef}
A Banach algebra $\A$ is called {\it weakly amenable\/} if ${\mathcal H}^1(\A,\A^\ast) = \{ 0 \}$.
\end{definition}
\begin{remark}
There is also the notion of a weakly amenable, locally compact group (\cite{CH}). This coincidence of terminology, however, is purely accidental.
\end{remark}
\par
In contrast to Theorem \ref{Barry1}, we have:
\begin{theorem}[\cite{Joh3}] \label{wamthm}
Let $G$ be a locally compact group. Then $L^1(G)$ is weakly amenable.
\end{theorem}
\par
For a particularly simple proof of this result, see \cite{DG}. For $M(G)$, things are strikingly different:
\begin{theorem}[\cite{DGH}] \label{DGHthm}
Let $G$ be a locally compact group. Then $M(G)$ is weakly amenable if and only if $G$ is discrete. In particular, $M(G)$ is amenable if and only if $G$ is discrete and amenable.
\end{theorem}
\par
Sometime after Kamowitz's pioneering paper, several mathematicians started to systematically develop a homological algebra with functional analytic overtones. Besides Johnson, who followed Hochschild's original approach, there were A.\ Guichardet 
(\cite{Gui}), whose point of view was homological rather than cohomological, and J.\ A.\ Taylor (\cite{Tay}) and --- most persistently --- A.\ Ya.\ Helemski\u{\i} and his Moscow school, whose approaches used projective or injective resolutions;
Helemski\u{\i}'s development of homological algebra for Banach and more general topological algebras is expounded in the monograph \cite{Hel2}.
\par
In homological algebra, the notions of projective, injective, and flat modules play a pivotal r\^ole. Each of these concepts tranlates into the functional analytic context. Helemski\u{\i} calls a Banach algebra $\A$ biprojective (respectively 
biflat) if it is a projective (respetively flat) Banach $\A$-bimodule over itself. We do not attempt to give the fairly technical definitions of a projective or a flat Banach $\A$-bimodule. Fortunately, there are equivalent, but more elementary 
characterizations of biprojectivity and biflatness, respectively.
\par
We use $\tensor_\gamma$ to denote the completed projective tensor product of Banach spaces. If $\A$ is a Banach algebra, then $\A \tensor_\gamma \A$ has a natural Banach $\A$-bimodule structure via
\[
  a \cdot (x \tensor y) := ax \tensor y \quad\text{and}\quad (x \tensor y) \cdot a =: x \tensor y a \qquad (a,x,y \in \A).
\]
This turns the multiplication operator
\[
  \Delta \!: \A \tensor_\gamma \A \to \A, \quad a \tensor b \mapsto ab
\]
into a homomorphism of Banach $\A$-bimodules.
\begin{definition} \label{bidef}
Let $\A$ be a Banach algebra. Then:
\begin{enumerate}
\item[(a)] $\A$ is called {\it biprojective\/} if and only if $\Delta$ has bounded right inverse which is an $\A$-bimodule homomorphism.
\item[(b)] $\A$ is called {\it biflat\/} if and only if $\Delta^\ast$ has bounded left inverse which is an $\A$-bimodule homomorphism.
\end{enumerate}
\end{definition}
\par
Clearly, biflatness is a property weaker than biprojectivity.
\par
The following theorem holds (\cite[Theorem 51]{Hel1}):
\begin{theorem}[A.\ Ya.\ Helemski\u{\i}] \label{biprojthm}
Let $G$ be a locally compact group. Then $L^1(G)$ is biprojective if and only if $G$ is compact.
\end{theorem}
\par
Again, a classical property of $G$ is equivalent to a cohomological property of $L^1(G)$. The question for which locally compact groups $G$ the Banach algebra $L^1(G)$ is biflat seems natural at the first glance. However,
any Banach algebra is amenable if and only if it is biflat and has a bounded approximate identity (\cite[Theorem Vii.2.20]{Hel2}). Since $L^1(G)$ has a bounded approximate identity for any $G$, this means that $L^1(G)$ is biflat precisely
when $G$ is amenable.
\par
Let $G$ be a locally compact group. A {\it unitary representation\/} of $G$ on a Hilbert space $\mathfrak H$ is a group homomorphism $\pi$ from $G$ into the unitary operators on $\mathfrak H$ which is continuous with respect to the given
topology on $G$ and the strong operator topology on ${\mathcal B}({\mathfrak H})$. A function 
\[
  G \to \comps, \quad x \mapsto \langle \pi(x) \xi, \eta \rangle
\]
with $\xi, \eta \in {\mathfrak H}$ is called a {\it coefficient function\/} of $\pi$.
\begin{example}
The {\it left regular representation\/} $\lambda$ of $G$ on $L^2(G)$ is given by
\[
  \lambda(x)\xi := L_{x^{-1}} \xi \qquad (x \in G, \, \xi \in L^2(G)).
\] 
\end{example}
\begin{definition}[\cite{Eym}] \label{Fourierdef}
Let $G$ be a locally compact group.
\begin{enumerate}
\item[(a)] The {\it Fourier algebra\/} $A(G)$ of $G$ is defined as
\[
  A(G) := \{ f \!: G \to \comps : \text{$f$ is a coefficient function of $\lambda$} \}.
\]
\item[(b)] The {\it Fourier--Stieltjes algebra\/} $B(G)$ of $G$ is defined as
\[
  B(G) := \{ f \!: G \to \comps : \text{$f$ is a coefficient function of a unitary representation of $G$} \}.
\]
\end{enumerate}
\end{definition}
\par
It is immediate that $A(G) \subset B(G)$, that $B(G)$ consists of bounded continuous functions, and that $A(G) \subset {\mathcal C}_0(G)$. However, it is not obvious that $A(G)$ and $B(G)$ are linear spaces, let alone algebras. Nevertheless, the 
following are true (\cite{Eym}):
\begin{itemize}
\item Let $\cstar(G)$ be the enveloping $\cstar$-algebra of the Banach $^\ast$-algebra $L^1(G)$. Then $B(G)$ can be canonically identified with $\cstar(G)^\ast$. This turns $B(G)$ into a commutative Banach algebra.
\item Let $\VN(G) := \lambda(G)^{\prime\prime}$ denote the {\it group von Neumann algebra\/} of $G$. Then $A(G)$ can be canonically identified with the unique predual of $\VN(G)$. This turns $A(G)$ into a commutative Banach algebra
whose character space is $G$.
\item $A(G)$ is a closed ideal in $B(G)$.
\end{itemize}
\par
If $G$ is an abelian group with dual group $\Gamma$, then the Fourier and Fourier--Stieltjes transform, respectively, yield isometric isomorphisms $A(G) \cong L^1(\Gamma)$ and $B(G) \cong M(\Gamma)$. Consequently, $A(G)$ is amenable for any
abelian locally compact group $G$. It doesn't require much extra effort to see that $A(G)$ is also amenable if $G$ has an abelian subgroup of finite index (\cite[Theorem 4.1]{LLW} and \cite[Theorem 2.2]{For}). On the other hand, every amenable 
Banach algebra has a bounded approximate identity, and hence Leptin's theorem (\cite{Lep}) implies that the ame\-na\-bi\-li\-ty of $A(G)$ forces $G$ to be amenable. Nevertheless, the tempting conjecture that $A(G)$ is amenable if and only if $G$ is 
amenable is false:
\begin{theorem}[\cite{Joh4}]
The Fourier algebra of\/ $\operatorname{SO}(3)$ is not amenable.
\end{theorem}
\par
This leaves the following intriguing open question:
\begin{question} \label{Q1}
Which are the locally compact groups $G$ for which $A(G)$ is amenable?
\end{question}
\par
The only groups $G$ for which $A(G)$ is known to be amenable are those with an abelian subgroup of finite index. It is a plausible conjecture that these are indeed the only ones. The corresponding question for weak amenability is open as well. 
B.\ E.\ Forrest has shown that $A(G)$ is weakly amenable whenever the principal component of $G$ is abelian (\cite[Theorem 2.4]{For}).
\par
One can, of course, ask the same question(s) for the Fourier--Stieltjes algebra:
\begin{question} \label{Q2}
Which are the locally compact groups $G$ for which $B(G)$ is amenable?
\end{question}
\par
Here, the natural conjecture is that those groups are precisely those with a compact, abelian subgroup of finite index. Since $A(G)$ is a complemented ideal in $B(G)$, the hereditary properties of amenability for Banach algebras 
(\cite[Theorem 2.3.7]{Run}) yield that $A(G)$ has to be amenable whenever $B(G)$ is. It is easy to see that, if the conjectured answer to Question \ref{Q1} is true, then so is the one to Question \ref{Q2}.
\par
Partial answers to both Question \ref{Q1} and Question \ref{Q2} can be found in \cite{LLW} and \cite{For}.
\section{and operator spaces}
Given any linear space $E$ and $n \in \posints$, we denote the $n \times n$-matrices with entries from $E$ by ${\mathbb M}_n(E)$; if $E = \comps$, we simply write ${\mathbb M}_n$. Clearly, formal matrix multiplication turns
${\mathbb M}_n(E)$ into an ${\mathbb M}_n$-bimodule. Identifying ${\mathbb M}_n$ with the bounded linear operators on $n$-dimensional Hilbert space, we equip ${\mathbb M}_n$ with a norm, which we denote by $| \cdot |$.
\begin{definition}
An {\it operator space\/} is a linear space $E$ with a complete norm $\| \cdot \|_n$ on ${\mathbb M}_n(E)$ for each $n \in \posints$ such that
\begin{equation} \tag{R 1} \label{R1}
  \left\| \begin{array}{c|c} x & 0 \\ \hline 0 & y \end{array} \right\|_{n+m} = \max \{ \| x \|_n, \| y \|_m \}
  \qquad (n,m \in \posints, \, x \in {\mathbb M}_n(E), \, y \in {\mathbb M}_m(E))
\end{equation}
and
\begin{equation} \tag{R 2} \label{R2}
  \| \alpha \cdot x \cdot \beta \|_n \leq | \alpha | \| x \|_n | \beta | \qquad (n \in \posints, \, x \in {\mathbb M}_n(E), \, \alpha, \beta \in {\mathbb M}_n).
\end{equation}
\end{definition}
\begin{example}
Let $\mathfrak H$ be a Hilbert space. The unique $\cstar$-norms on ${\mathbb M}_n({\mathcal B}({\mathfrak H})) \cong {\mathcal B}({\mathfrak H}^n)$ turn ${\mathcal B}({\mathfrak H})$ and any of its subspaces into operator spaces.
\end{example}
\par
Given two linear spaces $E$ and $F$, a linear map $T \!: E \to F$, and $n \in \posints$, we define the the {\it $n$-th amplification\/} $T^{(n)} \!: {\mathbb M}_n(E) \to {\mathbb M}_n(F)$ by applying $T$ to each matrix entry.
\begin{definition} \label{cbdef}
Let $E$ and $F$ be operator spaces, and let $T \in {\mathcal B}(E,F)$. Then:
\begin{enumerate}
\item[(a)] $T$ is {\it completely bounded\/} if
\[
  \| T \|_{\mathrm{cb}} := \sup_{n \in \posints} \| T^{(n)} \|_{{\mathcal B}({\mathbb M}_n(E),{\mathbb M}_n(F))} < \infty.
\]
\item[(b)] $T$ is a {\it complete contraction\/} if $\| T \|_{\mathrm{cb}} \leq 1$.
\item[(c)] $T$ is a {\it complete isometry\/} if $T^{(n)}$ is an isometry for each $n \in \posints$.
\end{enumerate}
\end{definition}
\par
The following theorem due to Z.-J.\ Ruan marks the beginning of abstract operator space theory:
\begin{theorem}[\cite{Rua1}] \label{Ruan1}
Let $E$ be an operator space. Then there is a Hilbert space $\mathfrak H$ and a complete isometry from $E$ into ${\mathcal B}({\mathfrak H})$.
\end{theorem}
\par
To appreciate Theorem \ref{Ruan1}, one should think of it as the operator space analogue of the elementary fact that every Banach space is isometrically isomorphic to a closed subspace of ${\mathcal C}(\Omega)$ for some compact
Hausdorff space $\Omega$. One could thus {\it define\/} a Banach space as a closed subspace of ${\mathcal C}(\Omega)$ some compact Hausdorff space $\Omega$. With this definition, however, even checking, e.g., that $\ell^1$ is a
Banach space or that quotients and dual spaces of Banach spaces are again Banach spaces is difficult if not imposssible.
\par
Since any $\cstar$-algebra can be represented on a Hilbert space, each Banach space $E$ can be isometrically embedded into ${\mathcal B}({\mathfrak H})$ for some Hilbert space $\mathfrak H$. For an operator space, it is not important that,
but {\it how\/} it sits inside ${\mathcal B}({\mathfrak H})$.
\par
There is one monograph devoted to the theory of operator spaces (\cite{ER}) as well as an online survey article (\cite{Wit}).
\par
The notions of complete boundedness as well as of complete contractivity can be defined for multilinear maps as well (\cite[p.\ 126]{ER}). Since this is somewhat more technical than Definition \ref{cbdef}, we won't give the details
here. As in the category of Banach spaces, there is a universal linearizer for the right, i.e.\ completely bounded, bilinear maps: the {\it projective operator space tensor product\/} (\cite[Section 7.1]{ER}), which we denote by $\Tensor$.
\begin{definition}
An operator space $\A$ which is also an algebra is called a {\it completely contractive Banach algebra\/} if multiplication on $\A$ is a complete (bilinear) contraction.
\end{definition}
\par
The universal property of $\Tensor$ (\cite[Proposition 7.1.2]{ER}) yields that, for a completely contractive Banach algebra $\A$, the multiplication induces a complete (linear) contraction $\Delta \!: \A \Tensor \A \to \A$.
\begin{example} \label{ccex}
\begin{enumerate}
\item For any Banach space $E$, there is an operator space $\max E$ such that, for any other operator space $F$, every $T \in {\mathcal B}(E,F)$ is completely bounded with $\| T \|_{\mathrm{cb}} = \| T \|$ (\cite[pp.\ 47--54]{ER}).
Given a Banach algebra $\A$, the operator space $\max \A$ is a completely contractive Banach algebra (\cite[p.\ 316]{ER}).
\item Any closed subalgebra of ${\mathcal B}({\mathfrak H})$ for some Hilbert space $\mathfrak H$ is a completely contractive Banach algebra.
\end{enumerate}
\end{example}
\par
To obtain more, more interesting, and --- in the context of abstract harmonic analysis --- more relevant examples, we require some more operator space theory.
\par
Given two operator spaces $E$ and $F$, let 
\[
  {\mathcal{CB}}(E,F) := \{ T \!: E \to F : \text{$T$ is completely bounded} \}.
\]
It is easy to check that ${\mathcal{CB}}(E,F)$ equipped with $\| \cdot \|_{\mathrm{cb}}$ is a Banach space. To define an operator space structure on $\mathcal{CB}(E,F)$, first note that ${\mathbb M}_n(F)$ is, for each $n \in \posints$, an
operator space in a canonical manner. The (purely algebraic) identification
\[
  {\mathbb M}_n({\mathcal{CB}}(E,F)) := {\mathcal{CB}}(E,{\mathbb M}_n(F)) \qquad (n \in \posints)
\] 
then yields norms $\| \cdot \|_n$ on the spaces ${\mathbb M}_n({\mathcal{CB}}(E,F))$ that satisfy (\ref{R1}) and (\ref{R2}), which is not hard to verify. 
\par
Since, for any operator space $E$, the Banach spaces $E^\ast$ and $\mathcal{CB}(E,\comps)$ are iso\-me\-tri\-cal\-ly isomorphic (\cite[Corollary 2.2.3]{ER}), this yields a canonical operator space structure on the dual Banach space of an operator space.
In partiuclar, the unique predual of a von Neumann algebra is an operator space in a canonical way.
\par
We shall see how this yields further examples of completely contractive Banach algebras.
\par
We denote the $W^\ast$-tensor product by $\bar{\tensor}$.
\begin{definition} \label{Hopf}
A {\it Hopf--von Neumann algebra\/} is a pair $(\M, \nabla)$, where $\M$ is a von Neumann algebra, and
$\nabla$ is a {\it co-multiplication\/}: a unital, injective, $w^\ast$-continuous $^\ast$-homomorphism $\nabla \!: \M \to \M \bar{\tensor} \M$ which is co-associative, i.e.\ the diagram
\[
  \begin{CD}
  \M @>\nabla>> \M \bar{\tensor} \M \\                                                                             
  @V{\nabla}VV                                          @VV{\nabla \tensor \id_\M}V                     \\
  \M \bar{\tensor} \M @>>{\id_\M \tensor \nabla}> \M \bar{\tensor} \M \bar{\tensor} \M
  \end{CD}
\]
commutes. 
\end{definition}
\begin{example} \label{wstarG}
Let $G$ be a locally compact group.
\begin{enumerate}
\item Define $\nabla \!: L^\infty(G) \to L^\infty(G \times G)$ by letting
\[
  (\nabla \phi)(xy) := \phi(xy) \qquad (\phi \in L^\infty(G), \, x,y \in G).
\]
Since $L^\infty(G) \bar{\tensor} L^\infty(G) \cong L^\infty(G \times G)$, this turns $L^\infty(G)$ into a Hopf--von Neumann algebra.
\item Let $\wstar(G)$ be the enveloping von Neumann algebra of $C^\ast(G)$. There is a canonical $w^\ast$-continuous homomorphism $\omega$ from $G$ into the unitaries of $\wstar(G)$
with the following universal property: For any unitary representation $\pi$ of $G$ on a Hilbert space, there is unique $w^\ast$-continuous $^\ast$-homomorphism
$\theta \!: \wstar(G) \to \pi(G)''$ such that $\pi = \theta \circ \omega$. Applying this universal property to the representation
\[
  G \to \wstar(G) \bar{\tensor} \wstar(G), \quad x \mapsto \omega(x) \tensor \omega(x)
\]
yields a co-multiplication $\nabla \!: \wstar(G) \to \wstar(G) \bar{\tensor} \wstar(G)$.
\end{enumerate}
\end{example}
\par
Given two von Neumann algebras $\M$ and $\N$ with preduals $\M_\ast$ and $\N_\ast$, their $\wstar$-tensor product $\M \bar{\tensor} \N$ also has a unique predual $(\M \bar{\tensor} \N)_\ast$. Operator space theory allows to identify this
predual in terms of $\M_\ast$ and $\N_\ast$ (\cite[Theorem 7.2.4]{ER}):
\[
  (\M \bar{\tensor} \N)_\ast \cong \M_\ast \Tensor \N_\ast.
\]
Since $\VN(G) \bar{\tensor} \VN(H) \cong \VN(G \times H)$ for any locally compact groups $G$ and $H$, this implies in particular that
\[
  A(G \times H) \cong A(G) \Tensor A(H).
\]
\par
Suppose now that $\M$ is a Hopf--von Neumann algebra with predual $\M_\ast$. The comultiplication $\nabla \!: \M \to \M \bar{\tensor} \M$ is $w^\ast$-continuous and thus the adjoint map of a complete contraction
$\nabla_\ast \!: \M_\ast \Tensor \M_\ast \to \M_\ast$. This turns $\M_\ast$ into a completely contractive Banach algebra.
In view of Example \ref{wstarG}, we have:
\begin{example} Let $G$ be a locally compact group.
\begin{enumerate}
\item The multiplication on $L^1(G)$ induced by $\nabla$ as in Example \ref{wstarG}.1 is just the usual convolution product. Hence, $L^1(G)$ is a completely contractive Banach algebra.
\item The multiplication on $B(G)$ induced by $\nabla$ as in Example \ref{wstarG}.2 is pointwise multiplication, so that $B(G)$ is a completely contractive Banach algebra. Since $A(G)$ is an ideal in $B(G)$ and 
since the operator space strucures $A(G)$ has as the predual of $\VN(G)$ and as a subspce of $B(G)$ coincide, $A(G)$ with its canonical operator space structure is also a completely contractive Banach algebra.
\end{enumerate}
\end{example}
\begin{remark}
Since $A(G)$ fails to be Arens regular for any non-discrete or infinite, amenable, locally compact group $G$ (\cite{For1}), it cannot be a subalgebra of the Arens regular Banach algebra $\mathcal{B}({\mathfrak H})$.
Hence, for those groups, $A(G)$ is not of the form described in Example \ref{ccex}.2.
\end{remark}
\par
We now return to homological algebra and its applications to abstract harmonic analysis.
\par
An operator bimodule over a completely contractive Banach algebra $\A$ is an operator space $E$ which is also an $\A$-bimodule such that the module actions of $\A$ on $E$ are completely bounded bilinear maps. One can then define
operator Hochschild cohomology groups $\mathcal{OH}^1(\A,E)$ by considering only completely bounded derivations (all inner derivations are automatically completely bounded). It is routine to check that the dual space of an operator $\A$-bimodule
is again an operator $\A$-bimodule, so that the following definition makes sense:
\begin{definition}[\cite{Rua2}] \label{opamdef}
A completely contractive Banach algebra $\A$ is called {\it operator amenable\/} if $\mathcal{OH}^1(\A,E^\ast) = \{ 0 \}$ for each operator $\A$-bimodule $E$.
\end{definition}
\par
The following result (\cite[Theorem 3.6]{Rua2}) shows that Definition \ref{opamdef} is indeed a good one:
\begin{theorem}[Z.-J.\ Ruan] \label{Ruan2}
Let $G$ be a locally compact group. Then $G$ is amenable if and only if $A(G)$ is operator amenable.
\end{theorem}
\begin{remark}
A Banach algebra $\A$ is amenable if and only if $\max \A$ is operator amenable (\cite[Proposition 16.1.5]{ER}). Since $L^1(G)$ is the predual of the abelian von Neumann algebra $L^\infty(G)$, the canonical operator space
structure on $L^1(G)$ is $\max L^1(G)$. Hence, Definition \ref{opamdef} yields no information on $L^1(G)$ beyond Theorem \ref{Barry1}.
\end{remark}
\par
The following is an open problem:
\begin{question}
Which are the locally compact groups $G$ for which $B(G)$ is operator amenable?
\end{question}
\par
With Theorem \ref{DGHthm} and the abelian case in mind, it is reasonable to conjecture that $B(G)$ is operator amenable if and only if $G$ is compact. One direction is obvious in the light of Theorem \ref{Ruan2}; a partial result
towards the converse is given in \cite{RSp}.
\par
Adding operator space overtones to Definition \ref{wamdef}, we define:
\begin{definition} \label{opwamdef}
A completely contractive Banach algebra $\A$ is called {\it operator weakly amenable\/} if $\mathcal{OH}^1(\A,\A^\ast) = \{ 0 \}$.
\end{definition}
\par
In analogy with Theorem \ref{wamthm}, we have:
\begin{theorem}[\cite{Spr}]
Let $G$ be a locally compact group. Then $A(G)$ is operator weakly amenable.
\end{theorem}
\par
One can translate Helemski\u{\i}'s homological algebra for Banach algebras relatively painlessly to the operator space setting: This is done to some extent in \cite{Ari} and \cite{Woo1}.
Of course, appropriate notions of projectivity and flatness play an important r\^ole in this operator space homological algebra. Operator biprojectivity and biflatness can be defined as in the classical setting,
and an analogue --- with $\Tensor$ instead of $\tensor_\gamma$ --- of the characterization used for Definition \ref{bidef} holds.
\par
The operator counterpart of Theorem \ref{biprojthm} was discovered, independently, by O.\ Yu.\ Aristov and P.\ J.\ Wood:
\begin{theorem}[\cite{Ari}, \cite{Woo2}] \label{OandP}
Let $G$ be a locally compact group. Then $G$ is discrete if and only if $A(G)$ is operator biprojective.
\end{theorem}
\par
As in the classical setting, both operator amenability and operator biprojectivity imply operator biflatness. Hence, Theorem \ref{OandP} immediately supplies examples of locally compact groups $G$ for which $A(G)$ is operator biflat,
but not operator amenable. A locally compact group is called a {\it $[\operatorname{SIN}]$-group\/} if $L^1(G)$ has a bounded approximate identity belonging to its center. By \cite[Corollary 4.5]{RX}, $A(G)$ is also operator biflat whenever
$G$ is a $[\operatorname{SIN}]$-group. It may be that $A(G)$ is operator biflat for every locally compact group $G$: this question is investigated in \cite{ARSp}.
\par
All these results suggest that in order to get a proper understanding of the Fourier algebra and of how its cohomological properties relate to the underlying group, one has to take its canonical operator space structure into account.

\begin{thebibliography}{Wit {\it et al.}}
%
\bibitem[Ari]{Ari} O.\ Yu.\ Aristov, Biprojective algebras and operator spaces. {\it J.\ Math.\ Sci.\ (New York)\/} {\bf 111\/} (2002), 3339--3386.
%
\bibitem[A--R--Sp]{ARSp} O.\ Yu.\ Aristov, V.\ Runde, and N.\ Spronk. Operator biflatness of the Fourier algebra. In preparation.
%
\bibitem[B--C--D]{BCD} W.\ G.\ Bade, P.\ C.\ Curtis, Jr., and  H.\ G.\ Dales, Amenability and weak amenability for Beurling and Lipschitz algebras. {\it Proc.\ London Math.\ Soc.\/}\ 
(3) {\bf 55\/} (1987), 359--377.
%
\bibitem[C--E]{CE} H.\ Cartan and S.\ Eilenberg, {\it Homological algebra\/}. Princeton University Press, Princeton, 1956.
%
\bibitem[C--H]{CH} M.\ Cowling and U.\ Haagerup, Completely bounded multipliers of the Fourier algebra of a simple Lie group of real rank one.
{\it Invent.\ Math.\/}\ {\bf 96\/} (1989), 507--549.
%
\bibitem[D--Gh--H]{DGH} H.\ G.\ Dales, F.\ Ghahramani, and A.\ Ya.\ Helemski\u{\i}, The amenability of measure
algebras. {\it J.\ London Math.\ Soc.\/}\ {\bf 66\/} (2002), 213--226.
%
\bibitem[Day]{Day} M.\ M.\ Day, Means on semigroups and groups. {\it Bull.\ Amer.\ Math.\ Soc.\/}\ {\bf 55\/} (1949), 1054--1055.
%
\bibitem[D--Gh]{DG} M.\ Despi\'c and F.\ Ghahramani, Weak amenability of group algebras of locally compact groups.
{\it Canad.\ Math.\ Bull.\/}\ {\bf 37\/} (1994), 165--167.
%
\bibitem[E--R]{ER} E.\ G.\ Effros and Z.-J.\ Ruan, {\it Operator spaces\/}. Clarendon Press, Oxford, 2000.
%
\bibitem[Eym]{Eym} P.\ Eymard, L'alg\`ebre de Fourier d'un groupe localement compact. {\it Bull.\ Soc.\ Math.\ France\/} {\bf 92\/} (1964), 181--236.
%
\bibitem[Fol]{Fol} G.\ B.\ Folland, {\it A course in abstract harmonic analysis\/}. CRC Press, Boca Raton, Florida, 1995.
%
\bibitem[For 1]{For1} B.\ E.\ Forrest, Arens regularity and discrete groups. {\it Pacific J.\ Math.\/}\ {\bf 151\/} (1991), 217--227. 
%
\bibitem[For 2]{For} B.\ E.\ Forrest, Amenability and weak amenability of the Fourier algebra. Preprint (2000).
%
\bibitem[Gre]{Gre} F.\ P.\ Greenleaf, {\it Invariant means on locally compact groups\/}. Van Nostrand, New York--Toronto--London, 1969.
%
\bibitem[Gui]{Gui} A.\ Guichardet, Sur l'homologie et la cohomologie des alg\`ebres de Banach. {\it C.\ R.\ Acad.\ Sci.\ Paris\/}, S\'er.\ A {\bf 262\/} (1966), 38--42.
%
\bibitem[Her]{Her} C.\ S.\ Herz, Harmonic synthesis for subgruops. {\it Ann.\ Inst.\ Fourier (Grenoble)\/} {\bf 23\/} (1973), 91--123.       
%
\bibitem[H--R]{HR} E.\ Hewitt and K.\ A.\ Ross, {\it Abstract harmonic analysis\/}, I and II. Springer Verlag, Berlin--Heideberg--New York, 1963 and 1970. 
%
\bibitem[Hel 1]{Hel1} A.\ Ya.\ Helemski\u{\i}, Flat Banach modules and amenable algebras. {\it Trans.\ Moscow Math.\ Soc.\/}\ {\bf 47\/} (1985), 199--224.
%
\bibitem[Hel 2]{Hel2} A.\ Ya.\ Helemski\u{\i}, {\it The homology of banach and topological algebras\/} (translated from
the Russian). Kluwer Academic Publishers, Dordrecht, 1989.
%
\bibitem[Hoch 1]{Hoch1} G.\ Hochschild, On the cohomology groups of an associative algebra. {\it Ann.\ of Math.\/}\ (2) {\bf 46\/} (1945), 58--67.
%
\bibitem[Hoch 2]{Hoch2} G.\ Hochschild, On the cohomology theory for associative algebras. {\it Ann.\ of Math.\/}\ (2) {\bf 47\/} (1946), 568--579.
%
\bibitem[Joh 1]{Joh1} B.\ E.\ Johnson, Cohomology in Banach algebras. {\it Mem.\ Amer.\ Math.\ Soc.\/}\ {\bf 127\/} (1972).
%
\bibitem[Joh 2]{Joh2} B.\ E.\ Johnson, Derivations from $L^1(G)$ into $L^1(G)$ and $L^\infty(G)$. In: J.\ P.\ Pier (ed.), {\it Harmonic analysis (Luxembourg, 1987)\/}, pp.\ 191--198. Lectures Notes in Mathematics {\bf 1359\/}.
Springer Verlag, Berlin--Heidelberg--New York, 1988.
%
\bibitem[Joh 3]{Joh3} B.\ E.\ Johnson, Weak amenability of group algebras. {\it Bull.\ London Math.\ Soc.\/}\ {\bf 23\/} (1991), 281--284.
%
\bibitem[Joh 4]{Joh4} {\sc B.\ E.\ Johnson}, Non-amenability of the Fourier algebra of a compact group. {\it J.\ London Math.\ Soc.\/}\ (2) {\bf 50\/} (1994), 361--374.
%
\bibitem[Kam]{Kam} H.\ Kamowitz, Cohomology groups of commutative Banach algebras. {\it Trans.\ Amer.\ Math.\ Soc.\/}\ {\bf 102\/} (1962), 352--372.
%
\bibitem[L--L--W]{LLW} A.\ T.-M.\ Lau, R.\ J.\ Loy, and G.\ A.\ Willis, Amenability of Banach and $\cstar$-algebras on locally compact groups. {\it Studia Math.\/}\ {\bf 119\/} (1996), 161--178.
%
\bibitem[Lep]{Lep} H.\ Leptin, Sur l'alg\`ebre de Fourier d'un groupe localement compact. {\it C.\ R.\ Acad.\ Sci.\ Paris\/}, S\'er.\ A {\bf 266\/} (1968), 1180--1182.
%
\bibitem[MacL]{MacL} S.\ MacLane, {\it Homology\/}. Springer Verlag, Berlin--Heidelberg--New York, 1995. 
%
\bibitem[Neu]{vonN} J.\ von Neumann, Zur allgemeinen Theorie des Ma{\ss}es. {\it Fund.\ Math.\/}\ {\bf 13\/} (1929), 73--116.
%
\bibitem[Pat]{Pat} A.\ L.\ T.\ Paterson, {\it Amenability\/}. American Mathematical Society, Providence, 1988.
%
\bibitem[Pie]{Pie} J.\ P.\ Pier, {\it Amenable locally compact groups\/}. Wiley-Interscience, New York, 1984.
%
\bibitem[R--St]{Rei} H.\ Reiter and J.\ D.\ Stegeman, {\it Classical harmonic analysis and locally compact groups\/}. Clarendon Press, Oxford, 2000.
%
\bibitem[Rua 1]{Rua1} Z.-J.\ Ruan, Subspaces of $\cstar$-algebras. {\it J.\ Funct.\ Anal.\/}\ {\bf 76\/} (1988), 217--230.
%
\bibitem[Rua 2]{Rua2} Z.-J.\ Ruan, The operator amenability of $A(G)$. {\it Amer.\ J.\ Math.\/}\ {\bf 117\/} (1995), 1449--1474.
%
\bibitem[R--X]{RX} Z.-J.\ Ruan and G.\ Xu, Splitting properties of operator bimodules and operator ame\-na\-bi\-li\-ty of Kac algebras. In: A.\ Gheondea,
R.\ N.\ Gologan and D.\ Timotin, {\it Operator theory, operator algebras, and related topics\/}, pp.\ 193--216. The Theta Foundation, Bucharest, 1997.
%
\bibitem[Run]{Run} V.\ Runde, {\it Lectures on amenability\/}. Lecture Notes in Mathematics {\bf 1774\/}. Springer Verlag, Berlin--Heidelberg--New York, 2002.
%
\bibitem[R--Sp]{RSp} V.\ Runde and N.\ Spronk, Operator amenability of Fourier-Stieltjes algebras. Preprint (2001). 
%
\bibitem[Spr]{Spr} N.\ Spronk, Operator weak amenability of the Fourier algebra. {\it Proc.\ Amer.\ Math.\ Soc.\/}\ {\bf 130\/} (2002), 3609--3617.      
%
\bibitem[Tay]{Tay} J.\ A.\ Taylor, Homology and cohomology for topological algebras. {\it Adv.\ in Math.\/}\ {\bf 9\/} (1970), 137--182.
%
\bibitem[Wei]{Wei} C.\ A.\ Weibel, {\it An introduction to homological algebra\/}. Cambridge University Press, Cambridge, 1994.   
%
\bibitem[Wit {\it et al.}]{Wit} G.\ Wittstock {\it et al.\/}, What are operator spaces? --- An online dictionary. URL: {\tt http://www.math.uni-sb.de/$^\sim$ag-wittstock/projekt2001.html}
(2001).   
%
\bibitem[Woo 1]{Woo1} P.\ J.\ Wood, {\it Homological algebra in operator spaces with applications to harmonic analysis\/}. Ph.D.\ thesis, University of
Waterloo, 1999.  
%
\bibitem[Woo 2]{Woo2} P.\ J.\ Wood, The operator biprojectivity of the Fourier algebra. {\it Canadian J.\ Math.\/}\ (to appear).                   
%
\end{thebibliography}
\end{document}